%% file: main.tex
\documentclass[opre,nonblindrev]{informs3} 

\DoubleSpacedXI 

\usepackage{endnotes}
\let\footnote=\endnote

\usepackage{natbib}
 \bibpunct[, ]{(}{)}{,}{a}{}{,}%
 \def\bibfont{\small}%
\TheoremsNumberedThrough     
\ECRepeatTheorems

\EquationsNumberedThrough    

\usepackage{mathrsfs}
\usepackage{steinmetz}
\usepackage{accents}
\usepackage{graphicx,scalerel}
\usepackage{amsmath}
\usepackage{hyperref}
\usepackage{dirtytalk}
\usepackage{array, booktabs,tabularx}
\usepackage{longtable}

\newcommand{\im}{\text{\textbf{i}}}

\newcommand\sbullet[1][.4]{\mathbin{\ThisStyle{\vcenter{\hbox{%
  \scalebox{#1}{$\SavedStyle\bullet$}}}}}%
}
\newcommand{\phasor}[1]{\underaccent{\sbullet}{#1}}

\newif\ifphasor
\phasorfalse

\begin{document}


\RUNAUTHOR{I. Aravena, D.K. Molzahn, S. Zhang et al.} 

\RUNTITLE{ARPA-E GO Competition Challenge 1 Overview} 

\TITLE{Recent Developments in Security-Constrained AC Optimal Power Flow: Overview of Challenge 1 in the ARPA-E Grid Optimization Competition} 

%



\ARTICLEAUTHORS{%
\AUTHOR{Ignacio Aravena\textsuperscript{1,$a$}, Daniel K. Molzahn\textsuperscript{2,$b$}, Shixuan Zhang\textsuperscript{2,$c$}}
\AFF{\textsuperscript{1}Lawrence Livermore National Laboratory, \textsuperscript{2}Georgia Institute of Technology, \\ \EMAIL{aravenasolis1@llnl.gov, molzahn@gatech.edu,  zhangshixuanus@gatech.edu}} 
\AUTHOR{Cosmin G. Petra\textsuperscript{1,$a$}}
\AFF{\textsuperscript{1}Lawrence Livermore National Laboratory, \textsuperscript{$a$}Team gollnlp \\ \EMAIL{petra1@llnl.gov}}
\AUTHOR{Frank E. Curtis\textsuperscript{3,b}, Shenyinying Tu\textsuperscript{4,b}, Andreas W{\"a}chter\textsuperscript{4,b}, Ermin Wei\textsuperscript{4,b}, Elizabeth Wong\textsuperscript{5,b}}
\AFF{\textsuperscript{3}Lehigh University, \textsuperscript{4}Northwestern University, \textsuperscript{5}University of California, San Diego, \textsuperscript{b}Team GO-SNIP \\ \EMAIL{fec309@lehigh.edu, shenyinyingtu2021@u.northwestern.edu, andreas.waechter@northwestern.edu, ermin.wei@northwestern.edu, elwong@ucsd.edu}} 
\AUTHOR{Amin Gholami\textsuperscript{2,c}, Kaizhao Sun\textsuperscript{2,c}, Xu Andy Sun\textsuperscript{6,c}}
\AFF{\textsuperscript{6}Masachusetts Institute of Technology, \textsuperscript{c}Team GMI-GO \\ \EMAIL{a.gholami@gatech.edu, ksun46@gatech.edu, sunx@mit.edu}} 
\AUTHOR{Stephen T. Elbert\textsuperscript{7,$d$}, Jesse T. Holzer\textsuperscript{7,$d$}, Arun Veeramany\textsuperscript{7,$d$},
\AFF{\textsuperscript{7}Pacific Northwest National Laboratory, \textsuperscript{$d$}Grid Optimization Competition Organization Team \\ \EMAIL{steve.elbert@pnnl.gov, jesse.holzer@pnnl.gov, arun.veeramany@pnnl.gov}}}
} 

\ABSTRACT{%
\textbf{\textsf{Abstract:}}
The optimal power flow problem is central to many tasks in the design and operation of electric power grids. This problem seeks the minimum cost operating point for an electric power grid while satisfying both engineering requirements and physical laws describing how power flows through the electric network. By additionally considering the possibility of component failures and using an accurate AC power flow model of the electric network, the security-constrained AC optimal power flow (SC-AC-OPF) problem is of paramount practical relevance. To assess recent progress in solution algorithms for SC-AC-OPF problems and spur new innovations, the U.S. Department of Energy's Advanced Research Projects Agency--Energy (ARPA-E) organized Challenge~1 of the Grid Optimization (GO) competition. This special issue includes papers authored by the top three teams in Challenge~1 of the GO Competition (Teams gollnlp, GO-SNIP, and GMI-GO). To introduce these papers and provide context about the competition, this paper describes the SC-AC-OPF problem formulation used in the competition, overviews historical developments and the state of the art in SC-AC-OPF algorithms, discusses the competition, and summarizes the algorithms used by these three teams. 
}%



\maketitle

%

\input{src/introduction}
\input{src/formulation}

\input{src/state-of-the-art}

\input{src/competition}
\input{src/approach_summary}
\input{src/conclusion}

\section*{Acknowledgements}
The authors greatly appreciate the efforts of all the competition organizers, especially Timothy Heidel, Kory Hedman, and Richard O'Neill, for their vision and work in designing and managing the GO competition.

%
%
%

\theendnotes


\bibliographystyle{informs2014}
\bibliography{ref/gosnip,ref/gollnlp,ref/gmigo,ref/temp,ref/dan,ref/go-platform}

\ECSwitch
\hspace*{-22pt}{\Large\textsf{\textbf{Electronic Companion to:}} \emph{Recent Developments in Security-Constrained AC Optimal Power Flow: Overview of Challenge~1 in the ARPA-E Grid Optimization Competition}}
\vspace{18pt}
\input{src/EC_phasor_example}
\input{src/EC_platform}

\end{document}

%% file: src/introduction.tex
\section{Introduction}

The Security-Constrained Alternating Current Optimal Power Flow (SC-AC-OPF) problem is central to almost all optimization models in electric power grids. 
This problem seeks optimal settings for power generation and grid control equipment so as to minimize operating cost while ensuring that the system can continue to be operated in the event of localized equipment failures.

The SC-AC-OPF problem is a large-scale, nonsmooth, nonconvex, and nonlinear optimization problem, which, in the simplified case without security constraints, has been shown to have multiple local minima \citep{Hiskens2001,bukhsh2013,molzahn2017space} and a corresponding feasibility problem that is strongly NP-hard \citep{bienstock2019strong}.  The nonlinearities of the SC-AC-OPF problem involve products and trigonometric functions, its constraints are nonconvex and nonsmooth, and the problem can have several million variables and constraints. Furthermore, SC-AC-OPF must be solved frequently every few minutes for each power system. These attributes make SC-AC-OPF an extremely difficult problem for which no effective exact method is known.

Rapidly emerging technologies like renewable generation, batteries, and electric vehicles motivate the development of new techniques for solving SC-AC-OPF problems. Further, improvements in SC-AC-OPF software have the potential to save \$6B to \$19B of electricity costs yearly in the US \citep{ch1foa}, even without taking into account the potential gains in reliability.

To advance the state-of-the-art in this high-risk, high-reward topic, the Advanced Research Projects Agency--Energy (ARPA-E) in the US Department of Energy (DOE) launched the Grid Optimization (GO) Competition -- Challenge 1, focusing on SC-AC-OPF, in 2018. This paper overviews the ARPA-E GO Competition Challenge 1 from the perspective of both the organizers and the top-three teams in the competition. We introduce the SC-AC-OPF problem from Challenge~1 to specialists in operations research (Section~\ref{sec:formulation}), present a historical review and discuss recent developments (Section~\ref{sec:state-of-the-art}), discuss the competition  (Section~\ref{sec:competition}), compare our teams' algorithms (Section~\ref{sec:approaches}), and provide concluding remarks (Section~\ref{sec:conclusion}).



%% file: src/formulation.tex

\section{Problem Formulation}
\label{sec:formulation}

The SC-AC-OPF problem seeks to minimize the production cost of generators to supply the load demands using an alternating current (AC) network. This section introduces the classical version of this problem to non-specialists with a basic understanding of direct current (DC) electrical circuits. We cover the most recurrent questions from researchers without a power engineering background during the competition, such as ``why are there angles?'' and ``why are there two types of power?''  We start by reviewing AC circuits and phasor notation. This section is similar in spirit to \citep{Frank2016, Bienstock2022}, but we take a deeper look at power in AC circuits. We then describe the elemental components of a power grid: buses, generators, loads, shunts, and branches. We finally present the SC-AC-OPF problem by bringing those components together.

\subsection{Alternating current circuits}
\label{subsec:phasors}


AC circuits operating in steady-state are similar to steady-state DC circuits in that their quantities of interest can be determined by solving algebraic, often linear, equations. The key difference between the two is that, while steady-state DC circuits can be algebraically described using only real numbers, AC circuits require complex numbers representing sinusoids, which we call \emph{phasors}. 
We next summarize these concepts.
For further details, the reader is referred to
Section \ref{sec:example_phasors} of the electronic companion material for this article as well as 
textbooks in circuit theory, such as \citep[Part 2]{Alexander2008}, \citep{Hayt2018}.

Steady-state voltages and currents in an AC power system are modeled via sinusoidal waveforms. We denote a voltage waveform as $v(t) = \sqrt{2} V \cos(\omega t)$, where $V$ is the root mean square (RMS) voltage and $\omega$ is the angular frequency. Power system models are derived using the fundamental equations for the resistors, inductors, and capacitors: $i_R(t) = v(t)/R$ (Ohm's Law), $i_L(t) = 1/L \cdot \int_{-\infty}^t v(\tau) d\tau$ (Faraday's Law), and $i_C(t) = C \frac{dv}{dt}(t)$ (Ampere-Maxwell's Law), where $R$, $L$, and $C$ denote resistance, inductance, and capacitance, respectively.

Applying these equations to the voltage waveform $v(t)$ reveals that resistors relate voltages and currents via scalar multiplication, i.e., $i_R(t) = \sqrt{2} \frac{V}{R} \cos(\omega t)$. Capacitors and inductors relate voltages and currents by both scalar multiplication and a $\pi/2$ phase shift (in opposite directions for inductors versus capacitors), i.e., $i_L(t) = \sqrt{2} \frac{V}{\omega L} \cos\left(\omega t - \frac{\pi}{2} \right)$ and $i_C(t) = \sqrt{2} V \, \omega C \, \cos\left(\omega t + \frac{\pi}{2} \right)$. Observe that the voltage and current waveforms all have the same frequency $\omega$.

These relationships motivate electrical engineers to represent AC circuit signals as elements of the complex plane, where multiplication and phase shift correspond to multiplication by complex numbers. Moreover, note that sinusoids can be orthogonally decomposed into components with phases $0$ and $\pi/2$. This decomposition naturally maps to the real and imaginary axes on the complex plane, enabling addition of sinusoids to take place via these orthogonal components. We refer to this representation of sinusoids as elements in the complex plane as \emph{phasors}.

Formally, we define the phasor transformation $\mathscr{F}_\omega: \mathcal{S}_\omega \rightarrow \mathbb{C}$, where $\mathcal{S}_\omega$ is the set sinusoidal functions with frequency $\omega$, as:
\begin{equation}
    \mathscr{F}_\omega\big(\sqrt{2} A \cos(\omega t + \phi)\big) = A \exp(\im\, \phi) \equiv A \phase{\phi},
    \label{eq:phasor_definition}
\end{equation}
where $\im = \sqrt{-1}$. For simplicity, we denote phasors using a dot below its corresponding symbol, e.g., $\mathscr{F}_\omega(v(t)) = \phasor{v}$. The phasor transformation is bijective with inverse $\mathscr{F}_\omega^{-1}(\phasor{x}) = \Re\big(\phasor{x} \cdot \sqrt{2} \exp(\im \omega t)\big)$.

The phasor representation of the voltage and current relationships for a generic circuit element $e$ (resistor, capacitor, or inductor) can be expressed as $\phasor{i}_e = Y_e \phasor{v}$, which is identical to Ohm's law, but with complex quantities, where $Y_e$ is called the \emph{admittance} of the element; $Y_R = 1/R$, $Y_L = -\im/(\omega L)$, and $Y_C = \im \omega C$. The addition of currents is performed along the real and imaginary axes in the same fashion as in DC circuits. We also note that parallel admittances are added directly to obtain the total admittance. Conversely, the total equivalent admittance for a series connection of admittances is the reciprocal of the sum of their reciprocals. The real and imaginary parts of an admittance $Y$ are the \emph{conductance} $G$ and \emph{susceptance} $B$, respectively, i.e., $G = \Re(Y)$ and $B = \Im(Y)$.

Power in AC circuits is also represented using complex quantities, though complex power is not a phasor as defined above. We first review instantaneous power. For an element with voltage $v(t) = \sqrt{2} V \cos(\omega t + \alpha)$ and current $i(t)=\sqrt{2} I \cos(\omega t + \beta)$, the instantaneous power is $p(t) = v(t) \cdot i(t)$ which, via the variable change $\tau = t + \alpha/\omega$ (time translation to align with voltage), yields:
\begin{equation}
    \begin{aligned}
        p(\tau) = v(\tau) \cdot i(\tau) \,=&~ \sqrt{2} V \cos(\omega \tau) \cdot \sqrt{2} I \cos(\omega \tau + \beta - \alpha) \\
         =&~ V I \cos(\alpha - \beta) \, \big( 1 + \cos(2 \omega \tau )\big) + V I \sin(\alpha - \beta) \, \cos(2 \omega \tau + \pi/2).
    \end{aligned}
    \label{eq:inst_power}
\end{equation}
Instantaneous power has three components, all proportional to the product of RMS voltage and current: (\textit{i}) one continuous component, which is also  proportional to the cosine of the phase difference between voltage and current; (\textit{ii}) one sinusoidal component, with double the frequency of the circuit proportional to the cosine of the phase difference; and (\textit{iii}) another sinusoidal component of double frequency, but proportional to the sine of the phase difference and with a $\pi/2$ phase shift with respect to the other sinusoidal component. Observe that only the first component can make actual work (spending energy over time), whereas the second is necessary for the first, and both become larger as the phase difference becomes closer to zero. The third component, on the other hand, cannot produce any work over time.
The first and second component of instantaneous power is only present for resistors (producing work), whereas inductors and capacitors only have the third component (storing and releasing energy from their magnetic and electric fields, respectively).

To obtain meaningful power quantities when working with phasors, we define complex power as $\phasor{s} = \phasor{v} \cdot \phasor{i}^*$, where $\phasor{i}^*$ denotes the complex conjugate of the current phasor, and hence we have:
\begin{equation}
    \phasor{s} = \phasor{v} \cdot \phasor{i}^* = V\phase{\alpha} \cdot I\phase{-\beta} = V I \phase{\alpha - \beta} = \underbrace{VI \cos(\alpha - \beta)}_{\textstyle p} +\, \im\, \underbrace{VI \sin(\alpha - \beta)}_{\textstyle q},
    \label{eq:complex_power}
\end{equation}
where $p$ is referred to as the \emph{active power}, representing the first and second components of instantaneous power, and $q$ is referred to as the \emph{reactive power}, representing the third component of instantaneous power. 

The analysis of large-scale power systems is built, mostly, out of the elements reviewed in this section. Aspects not reviewed here, which are also relevant in the study of power systems, are the treatment of \emph{transformers}, \emph{three-phase circuits}, the \emph{per-unit system}, and \emph{symmetrical components}. For the purposes of this paper, it suffices to say that the per-unit system and symmetrical components can be used to obtain a single-phase equivalent power system model. When applied to this model, the material reviewed above provides the basis for understanding the SC-AC-OPF problem considered in Challenge~1 of the GO Competition. For further details on these topics, the reader is referred to power systems textbooks, such as \citep{Gonen2016, Glover2016}.

\subsection{Power grid components}

This section reviews the components of the power system in the OPF problem, their power and voltage variables, and the constraints on those variables. Each of these components can have two modes of operation, one for base-case (normal) conditions and another for post-contingency (emergency) conditions, as well as constraints indicating how to transition from base case to post-contingency mode. A contingency event, in this context, corresponds to the abrupt disconnection of a generator, line, or transformer, leaving the rest of the system to adjust to new operating conditions in order to avoid a cascading failure or blackout. The base to post-contingency transition constraints presented in this section assume that no decision can be made in between base-case and post-contingency states, but rather the transition is dictated by automatic controller actions based on local observations. This is commonly referred to as a \emph{preventive} model in the literature.\endnote{Whereas preventive models accurately represent the existing operation paradigm and capabilities of power grids, so-called \emph{corrective} models are also prevalent in the literature. Corrective models assume that decisions can be made and coordinated between base case and post-contingency states, and thus problems can be \emph{corrected} as supposed to only \emph{prevented}. 
}

We denote base-case parameters and variables with a $0$ subscript and contingency parameters and variables with a $k \in \mathcal{K}$ subscript, where $\mathcal{K} \not\owns 0$ is the set of contingencies.  For notational convenience, we let $\mathcal{K}_0 := \mathcal{K} \cup \{0\}$ denote the set of all system conditions.

From this section onward, we use the same notation as in ARPA-E's GO Competition Challenge~1 \citep{GOChallenge1_formulation}, except that we use upper case for parameters and lower case for variables.

\subsubsection{Buses}

Buses correspond to the nodes in the electrical network where generators, loads, branches, and shunts connect to the rest of the grid. Each bus $n$ in the set of buses $\mathcal{N}$ has voltage phasors $v_{n,0} \phase{\theta_{n,0}}$ in the base case and $v_{n,k} \phase{\theta_{n,k}}$ in each post-contingency condition $k \in \mathcal{K}$, with (RMS) voltage magnitude limits:
\begin{gather}
    \underline{V}_n \leq v_{n,0} \leq \overline{V}_n ~~ \forall n \in \mathcal{N}, \label{eq:v_bnds_base} \\
    \underline{V}^E_n \leq v_{n,k} \leq \overline{V}^E_n ~~ \forall n \in \mathcal{N}, k \in \mathcal{K}, \label{eq:v_bnds_contingency}
\end{gather}
where $\underline{V}^E_n \leq \underline{V}_n \leq \overline{V}_n \leq \overline{V}^E_n$ for all $n \in \mathcal{N}$ are lower and upper voltage bounds. Since angles only make sense relative to one another, a reference bus $n^{\text{ref}} \in \mathcal{N}$ is chosen to have $\theta_{n^\text{ref},k} = 0$ for all $k \in \mathcal{K}_0$.


\subsubsection{Generators}

Generators transform primary energy sources into electrical power that is injected into the power system, that is, they act as voltage sources from Section \ref{subsec:phasors}. Each generator~$g$ in the set of generators $\mathcal{G}$ can inject active power $p_{g,0}$ (base case), $p_{g,k}$ (post-contingency) and reactive power $q_{g,0}$ (base case), $q_{g,k}$ (post-contingency) into the network within certain limits:
\begin{gather}
    \underline{P}_g \leq p_{g,k} \leq \overline{P}_g,~~ \underline{Q}_g \leq q_{g,k} \leq \overline{Q}_g ~~ \forall g \in \mathcal{G}, k \in \mathcal{K}_0 \label{eq:pq_bnds} 
\end{gather}
Let $n(g)$ denote the bus where the generator is connected. In the base case, both active and reactive injections can be decided freely, within the above limits, irrespective of the voltage magnitude at the generator's bus, $v_{n(g),0}$. In post-contingency state $k \in \mathcal{K}$, however, if the generator is in the set of generators that remain online $\mathcal{G}(k) \subseteq \mathcal{G}$, the generator's behavior is constrained by two controllers.

First, the voltage control system attempts to maintain the base-case voltage magnitude $v_{n(g),0}$ in the post-contingency state by injecting or absorbing reactive power; voltage magnitude increases with reactive power injection. The voltage control system only allows the post-contingency voltage $v_{n(g),k}$ to deviate from $v_{n(g),0}$ after depleting the generator's reactive power capability, fixing $q_{g,k}$ to its upper limit if voltage decreases post-contingency and to its lower limit if voltage increases post-contingency \citep{stott1974}. This behavior is modeled via the following complementarity constraints:
\begin{gather}
    \nu^+_{n(g),k} - \nu^-_{n(g),k} = v_{n(g),k} - v_{n(g),0} \label{eq:v_regulator_aux} \\
	0 \leq \nu^-_{n(g),k} \perp \overline{Q}_g - q_{g,k} \geq 0 ~~ \forall g \in \mathcal{G}(k), k \in \mathcal{K} \label{eq:v_regulator_up} \\
	0 \leq \nu^+_{n(g),k} \perp q_{g,k} - \underline{Q}_{g} \geq 0 ~~ \forall g \in \mathcal{G}(k), k \in \mathcal{K} \label{eq:v_regulator_down}
\end{gather}
where $\nu^+_{n(g),k}$ and $\nu^-_{n(g),k}$ are the voltage increase and decrease, respectively, 
in post-contingency $k$.

Second, the droop control system, also referred to as the automatic generation control system, adjusts the post-contingency active power $p_{g,k}$ following the post-contingency deviation in system frequency $\delta_k$ for $k \in \mathcal{K}$. This system increases generation if there is a frequency decrease ($\delta_k < 0$, indicating a lack of generation) and decreases it in case of a frequency increase ($\delta_k > 0$, indicating a generation excess) \citep{jaleeli1992}. Each generator $g$ may have a different droop slope $A_g \geq 0$, indicating how much to increase or decrease production per each unit of frequency deviation. The droop control saturates at the active power limits, as a generator can only operate within them. This behavior can be modeled as the following complementarity constraints:
\begin{gather}
    \rho^+_{g,k} - \rho^-_{g,k} = p_{g,k} - (p_{g,0} + A_{g} \delta_{k}) \label{eq:droop_control_auxiliary} \\
	0 \leq \rho^-_{g,k} \perp \overline{P}_g - p_{g,k} \geq 0 ~~ \forall g \in \mathcal{G}(k), k \in \mathcal{K} \label{eq:droop_control_up_saturation} \\
	0 \leq \rho^+_{g,k} \perp p_{g,k} - \underline{P}_g \geq 0 ~~ \forall g \in \mathcal{G}(k), k \in \mathcal{K} \label{eq:droop_control_down_saturation}
\end{gather}
where $\rho^+_{g,k}$ and $\rho^-_{g,k}$ are the upward and downward deviations from the linear response of generator~$g$ in post-contingency $k$ due to saturation at minimum and maximum power output, respectively. If generator $g \in \mathcal{G}$ is not online in post-contingency $k \in \mathcal{G}$, i.e., $g \notin \mathcal{G}(k)$, then it cannot provide power:
\begin{equation}
    p_{g,k} = q_{g,k} = 0 ~~ \forall g \in \mathcal{G}\setminus\mathcal{G}(k), k \in \mathcal{K} \label{eq:failed_generators}
\end{equation}

\subsubsection{Loads}

Loads are modeled as constant withdrawals of active and reactive power from the system. The load at bus $n \in \mathcal{N}$ is denoted as $P^L_n$ for active power and $Q^L_n$ for reactive power.

\subsubsection{Shunts}

In the context of OPF, shunts correspond to variable inductors or capacitors---or a mixture of the two---connecting between a bus and the ground, thereby ``deviating'' current from its main path from generator to consumers.\endnote{Shunts can also be connected between different phases of a bus in a three-phase power system. These shunts can be equivalently represented by a shunt connected between the bus and the ground, or a neutral point, using the $\Delta$-Y transformation.} Shunts are used to withdraw (inductor) or inject (capacitor) reactive power to the grid, helping with reactive power and voltage control.

The total susceptance of shunts at bus $n \in \mathcal{N}$ is denoted $b_{n,k}$ for $k \in \mathcal{K}_0$ and is bounded as follows:
\begin{gather}
    \underline{B}_n \leq b_{n, k} \leq \overline{B}_n ~~ \forall n \in \mathcal{N}, k \in \mathcal{K}_0 \label{eq:b_bnds} 
\end{gather}

\subsubsection{Branches}


Branches correspond to the paths traversed by electricity from generators to loads. This paper covers AC transmission lines, as they are the most prevalent branch in existing power systems. Other types of branches, not covered here, include transformers, phase-shifters, and DC transmission lines; the reader is referred to \citep{Gonen2016} for further information.

\begin{figure}
    \centering
    \includegraphics{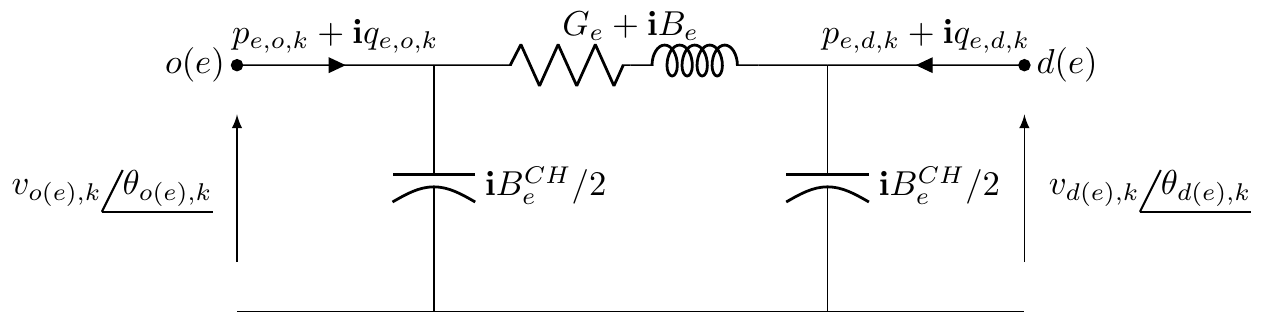}
    \caption{Equivalent circuit for transmission line $e \in \mathcal{E}(k)$ ($k \in \mathcal{K}_0$), from bus $o(e) \in \mathcal{N}$ to bus $d(e) \in \mathcal{N}$. The circuit uses one resistor and one inductor, with combined admittance $G_e + \im B_e$, to model the thermal losses and magnetic fields around the conductors. The circuit also uses two capacitors, with admittance $\im B_e^{CH}/2$ each, to model the effect of electric fields between the conductor and the ground. Complex power is defined as entering the line at both $o(e)$ and $d(e)$.}
    \label{fig:line_circuit}
\end{figure}

Each transmission line $e$ in the set of lines $\mathcal{E}$ has two terminals,  origin $o$ and destination $e$, which are connected to buses $o(e) \in \mathcal{N}$ and $d(e) \in \mathcal{N}$, respectively.  Each line $e$ is modeled using the circuit presented in Fig. \ref{fig:line_circuit}.
The complex power from bus $o(e)$ entering line $e$ through its $o$ terminal in system condition $k \in \mathcal{K}_0$ can be expressed as
\begin{equation*}
    \begin{aligned}\small
        p_{e,o,k} + \im q_{e,o,k} =&~ v_{o(e),k}\phase{\theta_{o(e),k}} \Bigg[ \bigg(\im \frac{B_e^{CH}}{2} \cdot  v_{o(e),k}\phase{\theta_{o(e),k}} \bigg)^* +  \Big((G_e + \im B_e) \big(v_{o(e),k}\phase{\theta_{o(e),k}} - v_{d(e),k}\phase{\theta_{d(e),0}}\big)\Big)^*\Bigg]
    \end{aligned}
\end{equation*}
which, when separated in real (active) and imaginary (reactive) parts, yields
\begin{align}
    &\begin{aligned}
        p_{e,o,k} =&\, G_e v_{o(e),k} - G_e \, \cos\big(\theta_{o(e),k} - \theta_{d(e),k}\big) \, v_{o(e),k} \, v_{d(e),k} - B_e \, \sin\big(\theta_{o(e),k} - \theta_{d(e),k}\big) \, v_{o(e),k} \, v_{d(e),k}
    \end{aligned} \label{eq:p_branch_o} \\
    &\begin{aligned}
        q_{e,o,k} =& -\big(B_e + B_e^{CH}/2\big) \, v_{o(e),k}^2 + B_e \, \cos\big(\theta_{o(e),k} - \theta_{d(e),k}\big) \, v_{o(e),k} \, v_{d(e),k}~- \\
        & \quad G_e \, \sin\big(\theta_{o(e),k} - \theta_{d(e),k}\big) \, v_{o(e),k} \, v_{d(e),k}
    \end{aligned} \label{eq:q_branch_o}
\end{align}
for all $e \in \mathcal{E}, k \in \mathcal{K}_0$. Analogous equations for the power into the $d$ terminal from the $d(e)$ bus are, for all $e \in \mathcal{E}, k \in \mathcal{K}_0$,
\begin{align}
    &\begin{aligned}
        p_{e,d,k} =&\, G_e v_{d(e),k} - G_e \, \cos\big(\theta_{d(e),k} - \theta_{o(e),k}\big) \, v_{d(e),k} \, v_{o(e),k} - B_e \, \sin\big(\theta_{d(e),k} - \theta_{o(e),k}\big) \, v_{d(e),k} \, v_{o(e),k} 
    \end{aligned} \label{eq:p_branch_d} \\
    &\begin{aligned}
        q_{e,d,k} =& -\big(B_e + B_e^{CH}/2\big) \, v_{d(e),k}^2 + B_e \, \cos\big(\theta_{d(e),k} - \theta_{o(e),k}\big) \, v_{d(e),k} \, v_{o(e),k} ~- \\
        & \quad G_e \, \sin\big(\theta_{d(e),k} - \theta_{o(e),k}\big) \, v_{d(e),k} \, v_{o(e),k}
    \end{aligned} \label{eq:q_branch_d}
\end{align}
Lines disconnected in post-contingency case $k \in \mathcal{K}$ cannot carry power, and thus
\begin{equation}
    p_{e,o,k} = p_{e,d,k} = q_{e,o,k} = q_{e,d,k} = 0 ~~ \forall e \in \mathcal{E} \setminus \mathcal{E}(k), k \in \mathcal{K} \label{eq:failed_lines}
\end{equation}

Equations \eqref{eq:p_branch_o}--\eqref{eq:q_branch_d} are derived using the polar representation of voltages and rectangular representation of power, yielding the so-called \emph{polar} formulation of the OPF problem. Alternatively, one can chose to represent voltage in rectangular coordinates, leading to the \emph{rectangular} formulation. Another choice is to use complex current instead of complex power, yielding the \emph{IV} formulation. Further changes of variables exist, yielding additional formulations of the OPF problem. We refer to \citep{molzahn_hiskens-fnt2019} for a details on alternative OPF formulations.

Each transmission line can admit a certain maximum current determined, primarily, by its thermal capacity. In order to limit the current magnitude using only complex power and voltage, observe that $\phasor{s} = \phasor{v} \cdot \phasor{i}^* \implies |\phasor{s}| = |\phasor{v}| \cdot |\phasor{i}|$, and thus one can write the current magnitude limit as:
\begin{gather}
    \sqrt{p_{e,o,0}^2 + q_{e,o,0}^2} \leq R_e v_{e,o,0} + \sigma_{e,o,0},~ \sqrt{p_{e,d,0}^2 + q_{e,d,0}^2} \leq R_e v_{e,d,0} + \sigma_{e,d,0} ~~ \forall e \in \mathcal{E} \label{eq:thermal_0}
\end{gather}
where $R_e > 0$ is the thermal rating of line $e$, and $\sigma_{e,o,0} \geq 0$ and $\sigma_{e,d,0} \geq 0$ are overload slack variables. For post-contingency situations, $R^E_e \geq R_e$ is defined as the emergency rating of line $e$ and
\begin{gather}
    \sqrt{p_{e,o,k}^2 + q_{e,o,k}^2} \leq R^E_e v_{e,o,k} + \sigma_{e,o,k}, ~
    \sqrt{p_{e,d,k}^2 + q_{e,d,k}^2} \leq R^E_e v_{e,d,k} + \sigma_{e,d,k} ~~ \forall e \in \mathcal{E}(k), k \in \mathcal{K} \label{eq:thermal_k}
\end{gather}
is enforced, where $\sigma_{e,o,k} \geq 0$ and $\sigma_{e,d,k} \geq 0$ are overload slack variables.

\subsection{SC-AC-OPF problem formulation}

Along with their individual constraints, all components of the power grid in conjunction must respect active and reactive power balance at each node $n \in \mathcal{N}$:
\begin{gather}
    \sum_{\substack{g \in \mathcal{G}(k): \\ n(g) = n}} p_{g,k} - P^L_n - 
		\sum_{\substack{e \in \mathcal{E}(k): \\ o(e) = n}} p_{e,o,k} - 
		\sum_{\substack{e \in \mathcal{E}(k): \\ d(e) = n}} p_{e,d,k} =
        \sigma^{P}_{n, k} ~~ \forall n \in \mathcal{N}, k \in \mathcal{K}_0
        \label{eq:p_balance} \\
	\sum_{\substack{g \in \mathcal{G}(k): \\ n(g) = n}} q_{g,k} - Q^L_n + b_{n,k} v_{n,k}^2 -
		\sum_{\substack{e \in \mathcal{E}(k): \\ o(e) = n}} q_{e,o,k} - 
		\sum_{\substack{e \in \mathcal{E}(k): \\ d(e) = n}} q_{e,d,k} = 
        \sigma^{Q}_{n, k} ~~ \forall n \in \mathcal{N}, k \in \mathcal{K}_0
        \label{eq:q_balance}
\end{gather}
where $\sigma_{n,k}^P$ and $\sigma_{n,k}^Q$, for all $k \in \mathcal{K}_0$, are slack balance variables, and $\mathcal{G}(0) = \mathcal{G}$ and $\mathcal{E}(0) = \mathcal{E}$.

Piecing all components together, the SC-AC-OPF model is:
\begin{equation}
    \begin{aligned}
        \min_{\substack{p,q,b,v,\theta \\ \delta,\nu,\rho,\sigma}} ~~& \sum_{g \in \mathcal{G}} C_g(p_{g,0}) + \sum_{e \in \mathcal{E}} (C_e(\sigma_{e,o,0}) + C_e(\sigma_{e,d,0})) + \sum_{n \in \mathcal{N}} \big(C^P_n\big(\sigma^P_{n,0}\big) + C^Q_n\big(\sigma^Q_{n, 0}\big)\big) ~+ \\
        & \frac{1}{|\mathcal{K}|} \Bigg(\sum_{e \in \mathcal{E}(k)} (C_e(\sigma_{e,o,k}) + C_e(\sigma_{e,d,k})) + \sum_{n \in \mathcal{N}} \big(C^P_n\big(\sigma^P_{n,k}\big) + C^Q_n\big(\sigma^Q_{n, k}\big)\big)\Bigg) \\
        \text{s.t.} ~~~& \text{Base case constraints, }k=0\text{: }\eqref{eq:v_bnds_base}, \eqref{eq:pq_bnds}, \eqref{eq:b_bnds}, \eqref{eq:p_branch_o}\text{--}\eqref{eq:q_branch_d}, \eqref{eq:thermal_0}, \eqref{eq:p_balance}, \eqref{eq:q_balance} \\
        & \text{Transition constraints, } k \in \mathcal{K}\text{: }\eqref{eq:v_regulator_up}\text{--}\eqref{eq:droop_control_down_saturation} \\
        & \begin{array}{ll}
            \text{Post-contingency constraints, } k \in \mathcal{K}\text{: } & \eqref{eq:v_bnds_contingency}, \eqref{eq:pq_bnds}, \eqref{eq:failed_generators}, \eqref{eq:b_bnds}, \eqref{eq:p_branch_o}\text{--}\eqref{eq:q_branch_d}, \eqref{eq:failed_lines}, \eqref{eq:thermal_k},  \eqref{eq:p_balance}, \eqref{eq:q_balance}
        \end{array}
    \end{aligned}  
    \label{eq:SC-AC-OPF}
\end{equation}
where $C_g$ is the generation cost function corresponding to fuel cost for thermal generators, opportunity costs for hydros, etc.; $C_e$ is a convex increasing penalization function for line overloading; and $C^P_n$ and $C^Q_n$ are convex penalization functions that are nonlinear and symmetric with respect to 0 
for power imbalance. In a nutshell, the SC-AC-OPF model \eqref{eq:SC-AC-OPF} seeks to minimize the operation cost for the base case (normal operation) while ensuring that there exists a feasible operating point to which the system would deviate in case of a contingency event.


%% file: src/state-of-the-art.tex

\section{History of Power Systems Optimization and Recent Developments}
\label{sec:state-of-the-art}

\subsection{Historical Review of Power System Optimization}
Power system optimization dates back to the earliest days of interconnected power grids. Since the early 1900s, system operators solved \emph{economic dispatch} problems to schedule the generators' power outputs in order to minimize system-wide operating costs~\citep{noakes1962}. Later formulations described by~\cite{stahl931} and~\cite{steinberg1943} included approximations of power losses in transmission lines. Operators solved economic dispatch problems using specialized computing machines (see~\cite{johnson1939}), potentially augmented with loss approximations from network analyzers, that is, analog computers that approximately modeled the transmission network~\citep{george1949}. With the advent of digital computing in the 1950s \citep{cohn2015}, economic dispatch formulations such as those considered by~\cite{early1953} and~\cite{squires1960} more accurately modeled power losses in transmission lines. 

Building on this foundation, J.L.~Carpentier at \'Electricit\'e de France first formulated what is now called the optimal power flow (OPF) problem~\citep{carpentier1962}. The OPF problem extends the economic dispatch problem to consider additional constraints such as line flow and voltage magnitude limits. The solving of OPF problems was motivated by the increasing interconnectivity of transmission systems~\citep{JulieCohn_2017}. Around the same time Carpentier formulated the OPF problem, erroneous programming of a single relaying device caused a widescale blackout of the Northeastern United States and Ontario, Canada in 1965. To avoid similar blackouts, T.~Dy-Liacco proposed an extension to OPF problems known as \emph{security constraints}. Security-constrained OPF (SC-OPF) problems ensure that an operating point will withstand all component failures within a specified set of \emph{contingencies} (often including the individual failures of each component, i.e., \emph{$N-1$ security}) without causing a cascade of further failures~\citep{dyliacco1967,stott1987}. (In this section, the terminology \emph{SC-OPF} refers to the general class of security-constrained optimal power flow problems that either consider an AC power flow model as in Challenge~1 of the GO Competition or instead focus on various approximated or relaxed power flow models, while \emph{SC-AC-OPF} refers specifically to SC-OPF problems that use an AC power flow model.) The importance of security was widely recognized; reflecting on another blackout, \cite{Carpentier1979a} noted ``\ldots \emph{the general failure of 19 December 1973 lasted only 3h but caused a loss of production for the country estimated to be at least the equivalent of 50 years' savings through economic dispatch} \ldots''~.

The first algorithms for SC-OPF problems were based on solving the Karush-Kuhn-Tucker conditions using diverse approximations for gradients and Hessians within first- and second-order methods~\citep{Cory1972, Alsac1974}. These approaches are effective for small systems, but do not scale to practically sized power grids. Subsequent advances during the 1980s came from the adoption of formal decomposition techniques in SC-OPF. In particular, Benders decomposition was employed to solve a simplified SC-OPF (linearized, active power only) with corrective actions (i.e.,  dispatch decisions in contingency subproblems) \citep{Pereira1985}. This work was later extended to the full nonlinear SC-AC-OPF problem using a generalized Benders decomposition technique \citep{Monticelli1987}, though this extension does not guarantee convergence due to the contingency subproblems' nonconvexity. During the 1980s, most industrial implementations were based on linear programming technology \citep{Alsac1990}, which was much more mature than nonlinear programming technology at that point. By the end of the 1980s, nonlinear programming techniques became competitive alternatives for SC-AC-OPF algorithms that showed promise for scaling up to industry requirements \citep{Papalexopoulos1989}. 


Despite continuing progress in the subsequent decades (see \cite{opf_survey,opf_litreview1993IandII,pscc2014survey} for OPF literature surveys), there are still several challenges that need to be addressed both in the formulation and solution of SC-OPF problems \citep{Momoh1997,stott2012,capitanescu2011survey,Capitanescu2016}. One key challenge is the nonconvex nature of the problem, which prevents its direct use in computing electricity prices (dual variables) and precludes global optimality guarantees. As will be discussed further in Section~\ref{subsec:simplified_opf}, operators have resorted to convex (usually linear) approximations of SC-AC-OPF problems to compute prices and solve globally the approximated problems while recovering feasible but suboptimal solutions for the original SC-AC-OPF problem in post-processing steps. Other challenges include the need to extend modeling capabilities (e.g., transient stability requirements \citep{Bruno2002,Yuan2003} and limits on the number of post-contingency corrective actions \citep{capitanescu2011survey,Phan2015}), scale to larger networks, and improve solution speed for the use of SC-AC-OPF algorithms in near-real-time operations.

Parallel computing has demonstrated promising results for addressing the aforementioned scalability and speed challenges. As the first parallel computing approach to SC-OPF, \cite{Rodrigues1994} decomposed the computations for a dual simplex method with constraint generation to solve a linear approximation of the problem. Subsequent work solved the nonlinear SC-AC-OPF problem in parallel by decomposing the factorizations in an interior point method \citep{Qiu2005,Petra_18_hiopdecomp}. Other recent parallel computing approaches, such as the work by \cite{Liu2013}, decompose the SC-AC-OPF problem at the formulation level by separating the problem into base-case and contingency subproblems. Subproblems are solved using off-the-shelf nonlinear programming packages such as Ipopt \citep{WaecBieg06}, and coordination among subproblems is achieved through first-order methods, similar to earlier work \citep{Pereira1985}.

With increasing computational capabilities, engineers have generalized SC-OPF formulations for applications beyond scheduling generator setpoints \citep{kirschen2018fundamentals}. Generalizations of SC-OPF problems are used to, e.g., optimize the statuses of switches that determine the network topology \citep{lyon2016}, choose the generators' on/off statuses in unit commitment problems \citep{padhy2004}, manage uncertainty from stochastic generators \citep{bienstock2014chance}, compute stability margins~\citep{avalos2009}, and plan system expansions \citep{mahdavi2019}.

\subsection{Approximating Optimal Power Flow Problems}\label{subsec:simplified_opf}
Many applications require solutions to SC-OPF problems within demanding time constraints, e.g., approximately 10 minutes for real-time applications. The computational difficulties of SC-AC-OPF problems have traditionally precluded direct solution within this time frame. Thus, operators often resort to power flow approximations that apply various assumptions based on characteristics of typical transmission systems to linearize the power flow equations~\citep{Stott1978,stott2009,molzahn_hiskens-fnt2019}. While these approximations provide advantages in computational speed and reliability, they may neglect important aspects of the power system physics (e.g., voltage magnitudes and reactive power) and  can introduce significant errors that impact both the achievable cost and the feasibility of the resulting solutions \citep{overbye2004hicss,purchala2005,coffrin2012,dvijotham_molzahn-cdc2016}. Nevertheless, operators have a long history of using power flow linearizations. For instance, the so-called ``DC power flow'' approximation was first used for fault analysis problems in the early 1900s \citep{wilson1916} and the New York Power Pool employed linearizations to optimize near-real-time generator dispatch in 1981~\citep{Elacqua1982}. Engineers have developed a wide variety of power flow approximations that model power losses in transmission lines, use proxy constraints to address more complex characteristics such as stability and voltage magnitude limits, and are tailored for particular applications.
\cite{molzahn_hiskens-fnt2019} provide a recent survey of power flow approximations.

Complementing these approximations, researchers have recently developed a wide range of power flow relaxations formulated as linear, second-order cone, and semidefinite programs \citep{low2014convex1, low2014convex2, molzahn_hiskens-fnt2019}. Relaxations enclose the nonconvex OPF feasible regions within a larger convex set and thus provide bounds on the globally optimal objective value. Moreover, whenever the solution to a convex relaxation is feasible for the original problem, the solution is guaranteed to be a global optimizer. Additionally, infeasibility of a convex relaxation certifies infeasibility of the nonconvex problem. At the same time, convex relaxations have a major practical disadvantage: they may have solutions that are infeasible for the original problem and recovering a feasible solution can be challenging \citep{venzke2020}.

In contrast to the outer enclosures used for convex relaxations, recently developed convex restrictions construct convex regions contained within an OPF problem's nonconvex feasible region. Using fixed-point theorems \citep{brouwer1911}, convex restrictions are constructed via a self-mapping set around a nominal operating point \citep{lee2019,cui2019}. Convex restrictions have been applied to compute feasible paths between operating points \citep{lee2019feasible} and solve robust AC-OPF problems that account for uncertain power injections \citep{lee2021robust_acopf}.

%% file: src/competition.tex


\section{Competition Description}
\label{sec:competition}
Recent developments in nonlinear optimization methods, alternative power flow modeling, and other technologies have raised the question of whether traditional approaches to SC-AC-OPF problems may be complemented or replaced by more sophisticated methods. Several optimization competitions have been organized to answer this question, notably including the Grid Optimization (GO) Competition run by the Advanced Research Projects Agency--Energy (ARPA-E) in the U.S. Department of Energy. This paper focuses on SC-AC-OPF problems in the form considered in Challenge 1 of the GO Competition as presented in Section~\ref{sec:formulation} \citep{GOChallenge1_formulation}. 

We note that subsequent ARPA-E competitions that consider different problem formulations are underway or forthcoming. Moreover, other power grid competitions have focused on applications of metaheuristics \citep{heuristic_competition} and machine learning \citep{marot2020learning,marot2021learning,l2rpn2020robustness,l2rpn2020adaptability,marot2021trust} as well as maintenance scheduling problems \citep{roadef,crognier2021}. More generally, competitions are frequently used to spur innovation in optimization, machine learning, and computer science. Existing platforms for conducting these competitions (e.g., \cite{kaggle}) do not provide computational resources on the scale needed by the GO Competition, especially timed runs on fixed hardware, thus necessitating customized hardware and software for this competition. The remainder of this section summarizes Challenge 1 of the ARPA-E GO Competition with a focus on aspects of the competition structure and computing platform that strongly influenced the solution algorithms.

With membership drawn from academia, national laboratories, and industry, 27 teams participated in the ARPA-E GO Challenge~1 competition. The teams had expertise in areas such as electric power systems, optimization, high-performance computing, and applied mathematics. Several teams included the developers of numerical optimization and modeling packages such as Ipopt~\citep{WaecBieg06}, SNOPT~\citep{gill2011}, Knitro~\citep{byrd2006}, BARON~\citep{sahinidis1996}, PowerModels.jl~\citep{coffrin2018}, and Gravity~\citep{hijazi2018}. 

The year-long competition involved three trial events before the final submission deadline, with increasingly difficult datasets released before and after each trial event for testing and evaluation. The teams' implementations produced outputs for 
the full solutions SC-AC-OPF problems within specified computational time requirements. These outputs were then scored by the competition organizers to determine each team's place within four different categories on the competition's leaderboard. We next summarize key aspects of the competition.\footnote{Further details on the competition structure are available at \url{https://gocompetition.energy.gov/challenges/challenge-1}.}

\subsection{Datasets}
\label{subsec:datasets}

%
%

%
%
%
%

The datasets provided in the competition consist of input files describing the network and scenarios for the load demands and generator costs. These datasets are formatted according to power system industry standards (namely, an extension of the data format used by the power system simulation tool PSS/E),
chosen for its familiarity throughout the field of power systems.\footnote{Details regarding the data format are available at \url{https://gocompetition.energy.gov/challenges/challenge-1/input-data-format}.} The trials and the final event used a mix of previously published and unseen networks, each of which included a variety of scenarios. The final event was conducted with 17 synthetic networks (each with 20 scenarios) and 3 actual industrial networks (each with 4 scenarios) for a total of 352 SC-AC-OPF problems. 

The industry datasets correspond to actual systems used in practical applications. As these datasets are confidential, the teams only received their overall scores for these datasets rather than detailed results. The synthetic datasets used in the competition were developed as part of the ARPA-E program ``Generating Realistic Information for the Development of Distribution and Transmission Algorithms'' (GRID DATA). Building on efforts dating back to at least the 1970s \citep{zaininger1977}, the GRID DATA program created \emph{realistic, but not real} datasets for use in the GO competition and other settings, i.e., datasets that are realistic enough to accurately benchmark algorithms, but free of confidential information and thus are publicly shareable.\footnote{For more information on the GRID DATA program, see \url{https://arpa-e.energy.gov/technologies/programs/grid-data}.} The methodologies for creating synthetic grids included ``top-down'' approaches that obfuscated information from actual industry datasets \citep{huang2018,soltan2019,fioretto2020,mak2020} and ``bottom-up'' approaches that built systems from scratch from geographic data such as population density and zoning information \citep{thiam2016,birchfield2017b}. Researchers have performed many studies assessing and validating the realism of these synthetic networks \citep{hines2010,birchfield2017a,birchfield2017b,pscc2018,li2018,sogol2021peci}. These and other datasets are available at two repositories, DR~POWER \citep{drpower} and BetterGrids \citep{nielsen2019}. Several datasets from the GO Competition are also in the PGLib-OPF library~\citep{pglib}.

The competition's datasets have networks ranging in size from 500 to 30,000 buses and include thousands of contingencies. The corresponding SC-AC-OPF problems thus have tens of thousands of variables and constraints in the base case alone. Considering contingencies yields problems with hundreds of millions of variables and constraints, thus posing significant computational challenges.

\subsection{Computing Platform and Architecture}
Teams participating in the GO Competition provided code implementing their algorithms to the competition organizers who then ran each team's code on the same computing platform. Each team accessed computing hardware, housed at Pacific Northwest National Laboratory (PNNL), with six 24-core computing nodes (two Intel Xeon E5-2670 v3 (Haswell) CPUs per node), each of which has 64 GB of memory, for a total of 144 cores.\footnote{See \url{https://gocompetition.energy.gov/evaluation-platform} for further details on the computing hardware.} 
This hardware was intended to represent the typical computing resources that an Independent System Operator (ISO) may have available to solve SC-AC-OPF problems in practical settings. The computing platform used the Centos operating system. Teams could access a number of optimization solvers (e.g., Gurobi, Cplex, Mosek, Ipopt, Knitro), programming languages (e.g., Python, Julia, Matlab), and modeling tools (e.g., GAMS, AMPL, CVX, JuMP).\footnote{For the final event, the teams submitted seven Python, four C++, and five each of Julia, executable binaries, and MATLAB codes.} Decomposition strategies and parallel computing (MPI between nodes and pthreads within each node) were key to efficiently using the computing hardware. 

\begin{figure}[htbp]
\centerline{\includegraphics[scale=0.33]{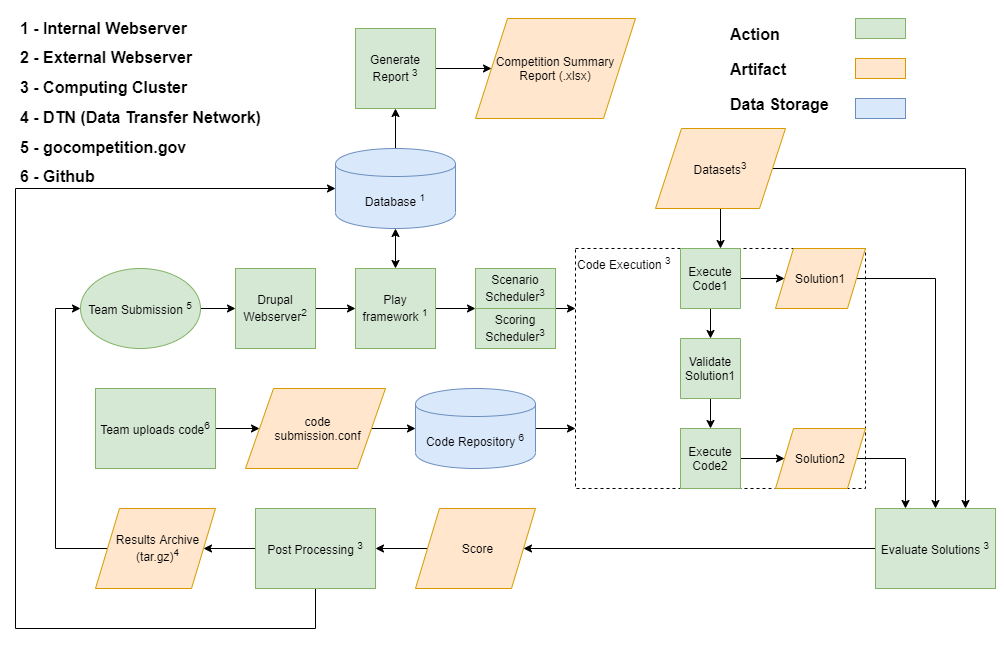}}
\caption{Workflow summary showing how teams' submissions were managed by the competition's organizers. Teams used a website to specify details like the programming language, location of their code, division of the competition, and datasets. The competition platform retrieved and executed each teams' code, evaluated the outputs to compute a score for each problem, and stored and reported the results. See Appendix~\ref{sec:platform-details} for further details.}
\label{fig:SysWorkflow}
\end{figure}

The full computing cluster at PNNL had over 500 nodes with up to 250 reserved for the exclusive use of the GO Competition. Challenge~1 accounted for over 200,000 computing tasks from 70,000 problems and used 16 million CPU-hours (wall clock) while generating over 200 petabytes of data. To manage these computing tasks, Fig.~\ref{fig:SysWorkflow} shows the competition organizer's workflow in executing the teams' codes and evaluating the results. Appendix~\ref{sec:platform-details} provides further details.

The organizers tracked hardware failures to automatically requeue and rerun affected submissions. Problems that timed out with no results were automatically rerun once. Computing runs that resulted in bad or infeasible results within the time limit were generally not rerun unless there was reason to suspect resource contention issues, e.g., similar runs on other problems were successful. Clusters of code failures also raised resource issue flags, the most common of which were network and disk access contention. Other challenges arose in ensuring consistency between runs executed on a team's computing hardware and the competition's platform. For instance, network traffic could lead to variability in run times. Additionally, differences in hardware and variations in software versions could affect the result of race conditions. In some cases, extensive troubleshooting between a team and the organizers was needed to resolve these and other issues.

\subsubsection{Required Outputs and Timing Requirements}
Each team's software had to provide values for the variables associated with both the base case and each contingency. However, the values for the base-case and post-contingency variables did not have to be computed simultaneously. Rather, the software was divided into two components: \emph{Code~1}, which determined the base-case solution, i.e., values of the variables associated with nominal operation, and \emph{Code~2}, which produced values for every contingency. More precisely, Code~1 observed the combined base-case and contingencies problem in order to optimize the base-case solution subject to the contingencies, but was only asked to provide a base-case solution.  Afterwards, Code~2 was run to produce values for post-contingency variables using the base-case solution that was determined by the run of Code~1.

The motivation for these two components of the software is that many industrial applications of SC-AC-OPF problems do not need explicit values for the post-contingency variables. Rather, the key requirement is obtaining base-case generator setpoints along with a strong degree of confidence that the base-case solution is secure against contingencies. This may not require explicit solutions for each contingency. However, the competition organizers needed post-contingency solutions to assess the security of the solution. With strict time limits for computing the base-case variables in Code~1 and permissive time limits for the post-contingency variables in Code~2, the competition assessed the practical relevance of the software with respect to both speed and security. 


The competition was also divided into two categories: real-time and offline optimization with strict 10-minute and 45-minute time limits, respectively, imposed on the execution of Code~1. These time requirements reflect the five- to fifteen-minute clearing times used in typical real-time electricity market operations as well as practical considerations related to assessing many scenarios in off-line applications. An additional time limit of two seconds per contingency was imposed on Code~2 regardless of category to prevent teams from submitting software with extremely long execution times while still maintaining tractability for computing the post-contingency variables.

%
%
%
%
%
%


\subsubsection{Solution Evaluation and Scoring Methods}
After applying a team's software to solve a set of problems, the organizers evaluated the solutions and scored the results. To evaluate a solution, hard constraints were checked for feasibility and the total objective was computed based on both the generation costs associated with power production and any penalty terms on the slack variables. These penalty terms ensured the feasibility of each team's outputs, despite any potential violations of power balance and branch flow limit constraints. The penalties were calculated via piecewise linear functions that significantly disincentivized large constraint violations. 

The choice to use penalty functions was motivated by the goal of evaluating and scoring as many solutions as possible. The competition organizers wanted to match a certain reality of industry practice. Specifically, if a solver fails for any reason, some kind of solution is still needed and the best solution available is used as far as possible. From this perspective, any ``infeasible'' solution is not really infeasible, but only worse than a feasible one to some degree. If a solution implies a power imbalance, then devices will adjust according to local control mechanisms in real time, e.g., generators will respond according to governor mechanisms and voltage-dependent loads will adjust their demands. Moreover, the data is also subject to error and uncertainty. However, a solution with no power imbalance is certainly preferable, so these penalties were especially severe.

Once the solutions for all the competitors over a set of problem instances were evaluated, the teams' algorithms were scored based on the evaluation results. The competition had two scoring methods. Based on similar benchmarking approaches for other classes of optimization problems \citep{mittelmann2022}, the first method computed the geometric mean of the scores for all scenarios associated with each network to compute a network score, then took the geometric mean of the network scores across all networks to compute an overall score for the team.
{\color{red}
%
%
}
The second method was based on the area under performance profiles that were constructed by comparing the ratios of teams' scores for each problem relative to the best score achieved for that problem by any team. The performance profile method was specifically intended to reward robust solver performance, at the possible expense of optimality, relative to the geometric mean method. Examination of the competition results showed that most teams did not change their algorithms in response to the different scoring methods. Moreover, analytical derivations showed that the two scoring methods were approximately equal to each other under a monotone transformation, so significantly different results could not have been expected.

A major challenge in scoring was how to handle the large differences in the scale of the objective values from one network to another. Generally, larger problems had larger optimal objective values, so simply adding up the objectives over all problem instances would incentivize good performance on the larger networks at the expense of performance on the smaller ones. Furthermore, for some problems, the gaps between the optimal objective and some easily computed relaxation objective or an objective corresponding to an easily computed heuristic solution might be much larger than the gaps for other problems. Implicitly, the objective of any optimization problem contains a constant term that is irrelevant to the quality of any given solution but affects both relative and absolute measures of optimality. The geometric mean scoring method was an attempt to deal with these issues in a way that ensured that each problem instance was meaningful to the overall score of a competitor.

Another major challenge in scoring involved solver outcomes that produced no solution, an unreadable or incorrectly formatted solution, an infeasible solution, or a very low quality solution. To address this challenge, the organizers computed a ``worst-case'' solution. The worst-case solution was created by keeping the optimization variables equal to their values in a pre-specified operating point, projecting to ensure feasibility of hard constraints, and accepting the resulting soft constraint violation penalties.\footnote{The detailed methodology for computing the worst-case score is available at \url{https://github.com/GOCompetition/WorstCase}.} In industry practice, a solver needs to have a highly reliable backup method that can be used if the main solver fails. The worst-case solution played that role in the competition. Solutions that either could not be scored otherwise or had poorer objective values than the worst-case score were assigned the worst-case score. Since these worst-case scores had a substantial negative impact on the overall results, teams needed to develop robust code that nearly always produced a correctly formatted solution of at least moderate quality within the specified time limit.


For prize awards, the competition had four divisions, corresponding to the two timing categories (10 minutes and 45 minutes) and the two scoring methods. Within each division, the top ten eligible teams were each awarded \$100,000, with prize eligibility based on a team's satisfaction of certain requirements including being led by a U.S. entity. An eligible team placing in the top ten in each of the four divisions would win \$400,000.


%% file: src/approach_summary.tex

%

\section{Summary of Approaches}
\label{sec:approaches}

This special issue includes three papers on algorithms for solving the SC-AC-OPF problem, as formulated in the GO Competition, from teams gollnlp, GO-SNIP, and GMI-GO. These teams ranked first, second, and third, respectively, in Challenge 1 out of the 27 participating teams. In this section, we briefly summarize the algorithms from each team and refer to the following references for further details~\citep{petra2021solving,curtis2021decomposition,gholami2022solving}.

The first team, gollnlp~\citep{petra2021solving}, uses nonconvex relaxations of complementarity constraints, two-stage decomposition with sparse approximation of recourse terms, and asynchronous parallelism, together with state-of-the-art nonlinear programming algorithms to compute SC-AC-OPF base case solutions that hedge against contingencies.
The gollnlp algorithm consistently provided high-quality feasible solutions to problems with different sizes and difficulty levels on both synthetic and industrial datasets.
The gollnlp algorithm produced the best-known solutions in $58\%$ of the cases. In the remaining cases, it attained average gaps (with respect to the best known solutions) of $0.15\%$ for real-time cases (with 10 minutes of execution time) and $0.21\%$ for offline cases (with 45 minutes of execution time). The gollnlp algorithm was also among the only three approaches that produced valid solutions for all cases.

The second team, GO-SNIP~\citep{curtis2021decomposition}, developed an algorithm based on an interior-point method that iteratively enforces constraints associated with additional contingencies.
This algorithm uses contingency screening and parallel processing techniques to identify quickly what seem to be the most important contingencies.
GO-SNIP also designed tailored heuristics for complementarity constraints and avoided certain degeneracies by modifying the problem formulation.
The GO-SNIP algorithm obtained the best-known solutions in 20\% of the cases, second-best-known solutions in over half of the cases, and top-ten solutions for almost all of the test cases. 

The third team, GMI-GO~\citep{gholami2022solving}, developed a two-level alternating direction method of multipliers (ADMM) algorithm for two-stage SC-AC-OPF problems with a convergence guarantee for the first time in the literature.
To handle the complementarity constraints, GMI-GO created a smoothing technique that allows interior-point solvers to be utilized.
To ensure robust performance, GMI-GO incorporated a contingency screening procedure and a massive parallel computation framework with safeguarding mechanisms.
The GMI-GO algorithm consistently outperformed the ARPA-E benchmark algorithm~\citep{coffrin2021arpa} in all instances and achieved the best-known solutions in about 5\% of the test cases.

The methods developed by these three teams share many similarities:
\begin{itemize}
    \item Due to the large sizes of the SC-AC-OPF problems, all three teams adopted decomposition algorithms or iterative contingency incorporation strategies, which are often effective in finding high-quality feasible solution within the given time limits.
    \item All three teams used contingency ranking and screening techniques in order to focus on those contingencies that seem to have a large influence on the base case decisions.
    \item As the original SC-AC-OPF problem formulation involves challenging nonlinear, nonconvex, and complementarity constraints, some approximation or relaxation techniques are considered in all three algorithms to simplify the SC-AC-OPF formulation.
    \item To better utilize the computational resources, each algorithm is built upon parallel and distributed computing methods adapted for the purpose of solving SC-AC-OPF problems.
    \item The implementation details also play an important role: for example, fallback mechanisms for failures in the solution process are used by each of the three teams.
    Moreover, the teams have all chosen C\texttt{++} as their programming language and used the interior-point solver Ipopt~\citep{WaecBieg06} with the HSL linear solver module~\citep{hsl2002}.
\end{itemize}

Despite these similarities, these algorithms are developed under distinctive assumptions and formulations, with different numerical performances on the test instances.
We therefore encourage readers to see each paper for further details regarding these algorithms \citep{petra2021solving,curtis2021decomposition,gholami2022solving}. We also note several resources published in other venues by GO Competition teams based on their work in Challenge 1. These include work by the sixth-place team Tartan Buffs \citep{bazrafshan2020computationally}, the eighth-place team Pearl Street Technologies \citep{jereminov2021equivalent}, and the tenth-place team ARPA-E Benchmark \citep{coffrin2018,coffrin2020code,coffrin2020obj}. The fifth-place team GravityX developed a modeling language, Gravity, to solve Challenge~1 \citep{hijazi2018}. 

%% file: src/conclusion.tex
\section{Conclusion}
\label{sec:conclusion}

Challenge 1 of the ARPA-E GO Competition compared the performance of a broad range of solution algorithms for SC-AC-OPF problems using state-of-the-art computing hardware and realistic, large-scale datasets. The results show that algorithms that combine power systems modeling techniques and physics-informed heuristics with advanced numerical optimization methods are capable of finding high-quality operating points satisfying the constraints of practical SC-AC-OPF problems within industrially relevant time frames. Thus, these results provide support for efforts to move nonlinear optimization techniques into industrial applications. This paper has given an overview of the Challenge~1 competition, including the problem formulation, competition structure, datasets, computing platform, etc., and summarized both the state of the art in SC-AC-OPF solution techniques and the algorithms used by the three top-performing teams in the competition. For further details of these algorithms, we refer to the three papers we authored for this special issue \citep{petra2021solving,curtis2021decomposition,gholami2022solving}. 

We conclude this introductory paper with comments on the advantages and limitations of conducting power systems research within the context of the competition and discuss directions for future work. The competition had the following significant advantages and research contributions.
%
\begin{itemize}
    \item \textbf{Algorithm development}: The competition encouraged teams to create new, highly effective algorithms for SC-AC-OPF problems that combine existing state-of-the-art techniques with novel approaches. These algorithms lay the foundation for addressing a wide range of challenging optimization tasks for future power grids.
    \item \textbf{Algorithm benchmarking}: The competition enabled benchmarking the performance of many SC-AC-OPF algorithms on common computing hardware, a consistent and practically relevant problem formulation, and the same scoring method. This fills a key gap in prior literature, as many earlier OPF algorithms consider incomparable problem formulations and are validated with inconsistent datasets, varying computing hardware, and limited comparisons to alternatives.
    \item \textbf{Building exposure}: The ARPA-E GO competition significantly raised the profile of SC-AC-OPF research and intensely focused the attention of researchers with a diversity of expertise on this topic. The competition thus provides the groundwork for future efforts in both subsequent GO Competition challenges and in other areas of power engineering research.
\end{itemize}

We also note several limiting aspects of performing research in the context of the competition.
\begin{itemize}
    \item \textbf{Prioritization of implementation over theory}: The GO Competition scoring method incentivizes teams to focus solely on reducing generation cost and avoiding constraint violations. These are important goals that are highly relevant to industrial applications. However, this neglects important complementary goals that focus more on theory, such as assessing the potential suboptimality of an operating point, proving infeasibility (i.e., certifying when the constraints in~\eqref{eq:SC-AC-OPF} cannot be satisfied with zero-valued slack variables), and bounding worst-case convergence rates. These and other goals have been the focus of intensive recent research~\citep{molzahn_hiskens-fnt2019, low2014convex1, low2014convex2}. Similarly, while the scoring method encourages developing heuristics that work well in many instances, there is little benefit for explaining and assessing the limitations to the performance of these heuristics within the context of the competition. Finally, since implementation quality was of utmost importance to avoid ``worst-case'' scores from code failures, the benchmark results from the competition may overlook alternative algorithms that could have performed better if they had more robust implementations or had better exploited the parallel computing resources.
    \item \textbf{Limited scope for off-line computations}: As discussed in Section~\ref{subsec:datasets}, the GO Competition evaluated the teams' implementations against previously unrevealed datasets, thus precluding the offline computations needed for many data-driven solution methodologies. However, in practice, many power grid parameters are relatively constant (e.g., line impedances) and system operators have good estimates for others (e.g., forecasts of load demands and renewable generation). Emerging machine learning algorithms perform off-line calculations using this data to inform on-line computations~\citep{duchesne2020}. Despite their significant promise, such approaches were not well suited for the GO Competition due to the competition's structure.
%
  In Challenge 2 of the GO Competition, a Monarch-of-the-Mountain (MoM) event---currently in progress---addresses
  this lack of practice with the data. MoM allows teams to use any hardware and any algorithm with no time limits,
  only evaluating solutions provided by the teams when they are ready.
  So far, with 62 scenarios having new solutions out of 84,
  only 2 show more than a 1\% improvement.
  These results suggest that the Challenge~2 teams did well without off-line computations.
    \item \textbf{Disincentives to disseminating in-progress research}: Many leading experts in power systems optimization were involved in the GO Competition as either members of various teams or in organizational roles. With competing teams and the organizers keeping their latest findings secret for the sake of fairness or competitive advantage, research for the competition differed substantially from the generally open environment that has historically been the norm in the power systems community. It is therefore essential that ongoing efforts such as the publications in this special issue as well as a special sessions at conferences and workshops continue to disseminate the knowledge developed in the GO competition.
\end{itemize}

Finally, we note several potential directions for further research. Challenges~2 and~3 of the GO Competition consider different problem formulations with features like multiple time periods, corrective controls, adjustable transformer tap settings, line switching, responsive load demands, and generator startup/shutdown characteristics. Research is also needed into other OPF problem formulations that are not considered in these challenges, such as uncertainties~\citep{pscc2022survey}, stability considerations~\citep{abhyankar2017,geng2017}, distributed solution algorithms~\citep{molzahn2017}, etc.

%% file: src/EC_phasor_example.tex
\section{Illustrative example of phasor transformation in AC circuits}
\label{sec:example_phasors}

\begin{figure}[h!]
    \centering
    \includegraphics{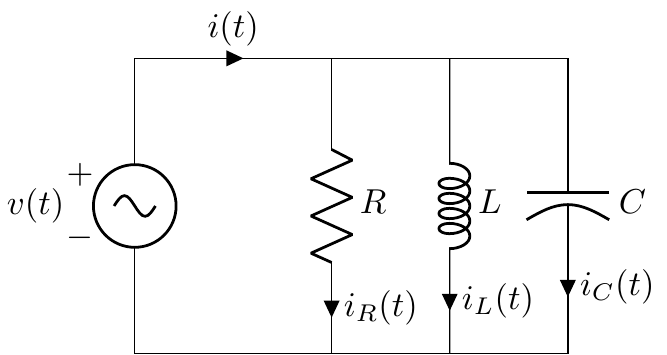}
    \caption{Example AC circuit with one voltage source, one resistor, one inductor, and one capacitor. $v(t)$ denotes the excitation voltage at the source (left), $i(t)$ the total current draw from the source, $i_R(t)$ the current passing by the resistor, which dissipates electrical power as heat, $i_L(t)$ the current by the inductor, which stores electrical energy in its magnetic field, and $i_C(t)$ the current by the capacitor, which stores electrical energy in its electrical field.}
    \label{fig:example_cirtuit}
\end{figure}

In this section, we demonstrate how phasors simplify computations in AC circuits, from differential equations in the time domain to algebraic equations in the complex plane, using the illustrative circuit of Fig. \ref{fig:example_cirtuit}.

\begin{figure}
    \centering
    \includegraphics{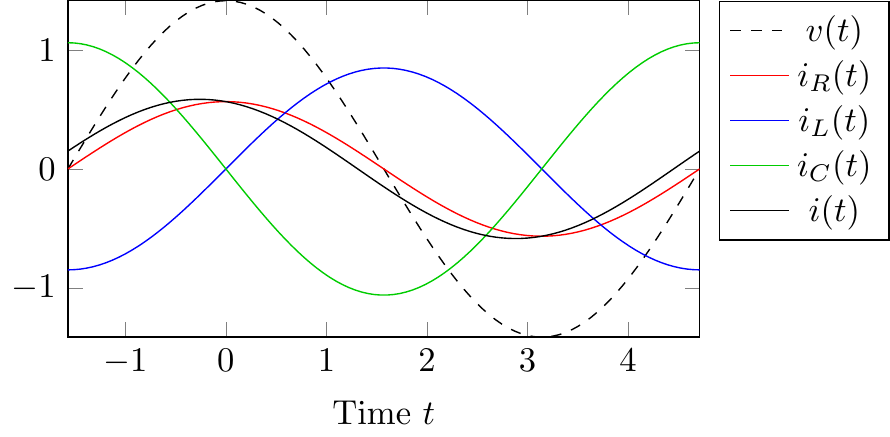}
    \caption{Oscillogram for voltages and currents in circuit of Fig. \ref{fig:example_cirtuit} with $R=2.5\Omega$, $L=1.3333$H, $C=0.75$F, and $v=\sqrt{2} \cos(t)$ ($V=1$ and $\omega = 1$).}
    \label{fig:oscillogram}
\end{figure}

The elemental circuit equations for passive elements (not producing power, i.e., resistor, inductor, and capacitor) imply: $i_R(t) = v(t)/R$ (Ohm's Law), $i_L(t) = 1/L \cdot \int_{-\infty}^t v(\tau) d\tau$ (Faraday's Law), $i_C(t) = C \frac{dv}{dt}(t)$ (Ampere-Maxwell's Law), and $i(t) = i_R(t) + i_L(t) + i_C(t)$ (Kirchhoff's Current Law). In steady-state, with an excitation $v(t) = \sqrt{2} V \cos(\omega t)$, where $V$ is the root mean square (RMS) voltage---which is easier to measure than the voltage amplitude---and $\omega$ is the angular frequency, these circuit equations simplify as follows:
\begin{align}
    i_R(t) =&~ \sqrt{2} \frac{V}{R} \cos(\omega t) \label{eq:ohm} \\
    i_L(t) =&~ \sqrt{2} \frac{V}{L} \cdot \int_{-\infty}^t \cos(\omega \tau) d\tau = \sqrt{2} \frac{V}{\omega L} \cos\Big(\omega t - \frac{\pi}{2} \Big) \label{eq:faraday} \\
    i_C(t) =&~ \sqrt{2} V C \frac{d(\cos(\omega t))}{dt} = \sqrt{2} V \, \omega C \, \cos\Big(\omega t + \frac{\pi}{2} \Big) \label{eq:ampere} \\
    i(t) =&~ i_R(t) + i_L(t) + i_C(t) = \sqrt{2} \frac{V}{R} \cos(\omega t) + \sqrt{2} V \Big(\omega C -  \frac{1}{\omega L}\Big) \, \cos\Big(\omega t + \frac{\pi}{2}\Big) \label{eq:kirchhoff}
\end{align}
We observe that circuit calculations for steady-state AC circuits correspond to multiplication by a scalar for resistors \eqref{eq:ohm}, or a multiplication by a scalar and a phase shift for inductors \eqref{eq:faraday} and capacitors \eqref{eq:ampere}, leaving the frequency $\omega$ unaltered. Further, addition is performed in orthogonal directions \eqref{eq:kirchhoff}: any sinusoidal can be decomposed into a part with phase $0$ and a part with phase $\pi/2$, and addition can be performed separately for phase $0$ and $\pi/2$ because they are orthogonal. Fig. \ref{fig:oscillogram} presents the waveforms of $i_R(t)$, $i_L(t)$, and $i_C(t)$ for particular values of the circuit parameters.

These observations motivate electrical engineers to represent AC circuit quantities as elements of the complex plane, where multiplication and phase shift correspond to multiplication by complex numbers, and where addition also takes place in orthogonal components. We map sinusoidal signals onto the complex plain using the phasor transformation $\mathscr{F}_\omega: \mathcal{S}_\omega \rightarrow \mathbb{C}$, where $\mathcal{S}_\omega$ is the set sinusoidal functions with frequency $\omega$, defined in \eqref{eq:phasor_definition}, and repeated here for convenience:
\begin{equation*}
    \mathscr{F}_\omega\big(\sqrt{2} A \cos(\omega t + \phi)\big) = A \exp(\im\, \phi) \equiv A \phase{\phi}
\end{equation*}

\begin{figure}
    \centering
    \includegraphics{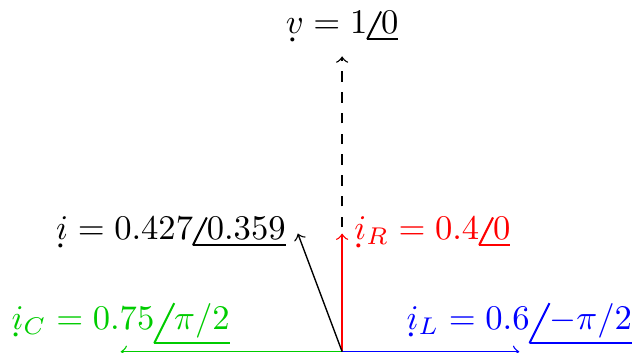}
    \caption{Phasor diagram for voltages and currents in circuit of Fig. \ref{fig:example_cirtuit}. Same electrical parameters as in Fig. \ref{fig:oscillogram}.}
    \label{fig:phasor_diagram}
\end{figure}

We apply the phasor transformation to the circuit equations \eqref{eq:ohm}--\eqref{eq:kirchhoff} to obtain the following compact, algebraic---as oppossed to differential---expressions:
\begin{align}
    \phasor{i}_R =&~ \frac{V}{R} \phase{0} \,=\, \frac{1}{R} \, \phasor{v} \label{eq:ohm_phasor} \\
    \phasor{i}_L =&~ \frac{V}{\omega L} \phase{-\frac{\pi}{2}} \,=\, \frac{-\im}{\omega L} \, \phasor{v} \label{eq:faraday_phasor} \\
    \phasor{i}_C =&~ \omega C V \phase{\frac{\pi}{2}} \,=\, \im \omega C \, \phasor{v} \label{eq:ampere_phasor} \\
    \phasor{i} =&~ \phasor{i}_R + \phasor{i}_L + \phasor{i}_C = \Bigg(\frac{1}{R} + \im\,\Big(\omega C -  \frac{1}{\omega L}\Big)\Bigg) \phasor{v} \label{eq:kirchoff_phasor}
\end{align}
The complex plane in Fig. \ref{fig:phasor_diagram} presents the phasors for the oscillogram from Fig.~\ref{fig:oscillogram}, where we can observe the scaling and rotation of the current with respect to the voltage for each passive element, as well as their aggregate.

For each passive element $e$, their phasor-equivalent equations \eqref{eq:ohm_phasor}--\eqref{eq:ampere_phasor} can be expressed as $\phasor{i}_e = Y_e \phasor{v}$, identical to Ohm's law but with complex quantities, where $Y_e$ is called the \emph{admittance} of the element; $Y_R = 1/R$, $Y_L = -\im/(\omega L)$, and $Y_C = \im \omega C$. Similarly, the addition of currents is performed along the real and imaginary axes in equation \eqref{eq:kirchoff_phasor}, in the same fashion that it would be done in DC circuits. We can also note that admittance in parallel can be added directly to obtain the total admittance $Y = 1/R + \im \, (\omega C - 1/(\omega L))$ (conversely, the total a of series admittance is the reciprocal of the sum of their reciprocals). We call the real and imaginary parts of an admittance $Y$ the \emph{conductance} $G$ and \emph{susceptance} $B$, respectively, i.e., $G = \Re(Y)$ and $B = \Im(Y)$.

As reviewed in Section \ref{subsec:phasors}, instantaneous power $p(t) = v(t) \cdot i(t)$ (equation \eqref{eq:inst_power}) can also be represented using complex quantities, and we defined complex power as $\phasor{s} = \phasor{v} \cdot \phasor{i}^*$ (equation \eqref{eq:complex_power}). We noted that the real component of complex power, active power, represents the continuous and double frequency part of instantaneous power proportional to the cosine of the phase difference between voltage $v(t)$ and current $i(t)$. The imaginary component of complex power, reactive power, represents the double frequency part of instantaneous power proportional to the sine of the phase difference between voltage $v(t)$ and current $i(t)$.

\begin{figure}
    \centering
    \includegraphics{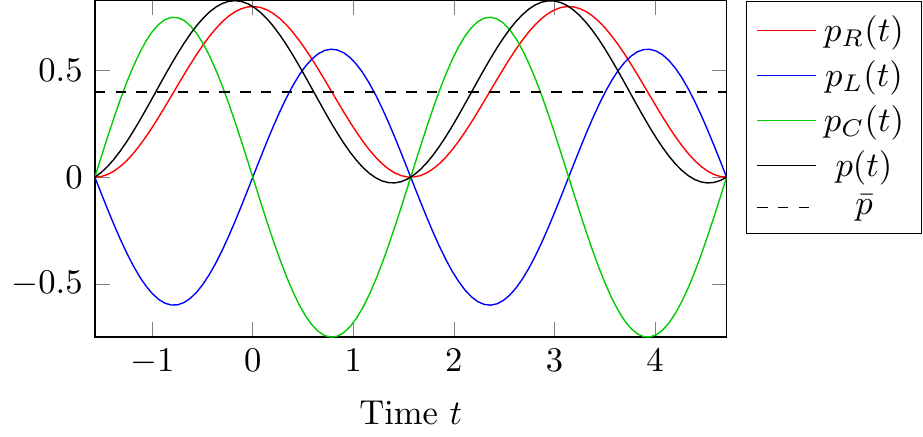}
    \caption{Oscillogram for instant power provided by the source, $p(t)$, and consumed by passive elements, $p_R(t)$, $p_L(t)$, and $p_C(t)$, in Fig. \ref{fig:example_cirtuit}. This oscillogram uses the same electrical parameters as in Fig. \ref{fig:oscillogram}.}
    \label{fig:oscillogram_power}
\end{figure}

\begin{figure}
    \centering
    \includegraphics{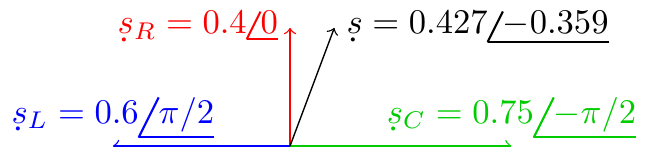}
    \caption{Complex power diagram for circuit in in Fig. \ref{fig:example_cirtuit} using the same electrical parameters as in Fig. \ref{fig:oscillogram}.}
    \label{fig:phasor_diagram_power}
\end{figure}

These relations can be visually appreciated in Figures \ref{fig:oscillogram_power} and \ref{fig:phasor_diagram_power}. Figure \ref{fig:oscillogram_power} presents the oscillogram for instantaneous power in the circuit of Fig. \ref{fig:example_cirtuit}. Figure \ref{fig:phasor_diagram_power} presents the complex power diagram for the circuit of Fig. \ref{fig:example_cirtuit}. For the resistor, instantaneous power $p_R(t)$ only has a continuous component of $0.4$ and double frequency sinusoidal component with amplitude $0.4$, both proportional to the magnitude of complex power $|\phasor{s}_R| = 0.4$ and the cosine of the angle of complex power $\angle\phasor{s}_R = 1$. Analog comparisons between instantaneous and complex power can be made for the inductor, capacitor, and the source. An additional, perhaps more important, observation from Fig. \ref{fig:phasor_diagram_power} is that that complex power also obeys balance as current does, that is $\phasor{s} = \phasor{s}_R + \phasor{s}_L + \phasor{s}_C$, which is a consequence of the definition of complex power.

%% file: src/EC_platform.tex
\section{Competition Platform Details}
\label{sec:platform-details}


This appendix summarizes the platform developed for Challenge 1 of the GO Competition. An internal web server is the central point of coordination among entrants, the computing cluster, database, code repository, and network file servers. 
The internal web server manages the running and evaluation of pending problems, receives progress updates while transitioning through the life cycle of a submission, and routes those messages to the entrant through a web interface. 

The web server validates and stores the details of each submission
and then triggers a chain of events. The team's code is checked out to the web server from GitHub using a predefined SSH key. 
If the submitted code does not yet exist on the login node of the computing cluster, it is copied from the web server and compiled if necessary. 
If this fails, logged information is sent to the team and the submission is marked as terminated. Otherwise, the code is then run.



Each submission includes a Code~1 which a 10-minute time limit for a real-time case and 45-minute limit for a planning case. Code~1 is expected to generate a \texttt{solution1.txt} file and optionally a \texttt{solution2.txt} file. If either Code~1 fails to produce the base case solution or produces an infeasible solution, the submission is terminated and assigned the worst-case score. Each submission may also include a Code 2 that is then run to generate a \texttt{solution2.txt} file using information from \texttt{solution1.txt}. Code~2 is not allowed to recreate or modify \texttt{solution1.txt} and has a time limit of two seconds per contingency. Once both solutions are available, the evaluation script is launched to generate a score. Once scored, post-processing archives results in a \texttt{tar.gz} file and stores scenario's metadata into a database. Should the team wish to review their score or obtain the results, they can do so through their submission page on the competition's web site.

Code 1 and Code 2 could be written in two different languages. The primary language was solicited from the submission webpage. 
The most common submissions were in the form of Linux binary executables, C/C++ (\href{https://gcc.gnu.org/}{gcc} versions 4.7/5.2/7.1/8.1), \href{www.mathworks.com}{MATLAB R2019a} with \href{https://matpower.org/}{MATPOWER v7.0b}, \href{https://julialang.org/}{Julia} (versions 0.64, 0.7. 1.1, 1.2) and \href{https://www.python.org/}{Python} (versions 2.7.13, 3.7.2). The entrants could choose from a broad range of compiler and solver versions, though latest versions of them were provided if there was no specific request. Entrants also had a wide range of choice of solver libraries compiled with various versions of gcc upon request, including \href{https://coin-or.github.io/Ipopt/index.html#Overview}{Ipopt} v3.12.13, \href{http://www.hsl.rl.ac.uk/ipopt/}{HSL for Ipopt} v2015.06, Python-based \href{http://www.pyomo.org/}{Pyomo} v5.6.1, Julia-based \href{http://www.juliaopt.org/JuMP.jl/dev/}{JuMP} (0.18, 0.19), \href{https://www.gurobi.com/}{Gurobi} v8.1, \href{https://www.ibm.com/analytics/cplex-optimizer}{IBM CPLEX} v12, MATLAB-based  \href{http://cvxr.com/cvx/doc/intro.html#what-is-cvx}{CVX} v2.1, \href{https://www.gams.com/}{GAMS} v27, \href{https://www.mosek.com/}{MOSEK} v7, \href{https://projects.coin-or.org/Bonmin}{Bonmin} v1, and \href{https://ampl.com/}{AMPL} v2019. In most cases,
commercial solver vendors sponsored a multi-use license for the length of the competition. Teams could also choose between Intel MPI (versions 2017 and 2018) and OpenMPI (versions 1.8.3 and 3.1) implementations with desired combinations of Intel and gcc compilers.